\documentclass[11pt,a4paper]{amsart}
\usepackage{amssymb,amsmath}
\usepackage{mathpazo}

\headsep=1cm \numberwithin{equation}{section}
\newtheorem{theorem}{Theorem}[section]
\newtheorem{lemma}[theorem]{Lemma}

\newtheorem{corollary}[theorem]{Corollary}

\theoremstyle{definition}
\newtheorem{definition}[theorem]{Definition}
\theoremstyle{remark}

\newtheorem{example}[theorem]{Example}

\newcommand{\Ass}{\operatorname{Ass}}
\newcommand{\im}{\operatorname{im}}
\newcommand{\grade}{\operatorname{grade}}

\newcommand{\Spec}{\operatorname{Spec}}

\newcommand{\ara}{\operatorname{ara}}

\newcommand{\cd}{\operatorname{cd}}

\newcommand{\Ht}{\operatorname{ht}}

\newcommand{\id}{\operatorname{id}}
\newcommand{\Gid}{\operatorname{Gid}}
\newcommand{\pd}{\operatorname{pd}}
\newcommand{\Gpd}{\operatorname{Gpd}}
\newcommand{\V}{\operatorname{V}}

\newcommand{\Ext}{\operatorname{Ext}}
\newcommand{\Supp}{\operatorname{Supp}}
\newcommand{\Tor}{\operatorname{Tor}}
\newcommand{\Hom}{\operatorname{Hom}}
\newcommand{\Att}{\operatorname{Att}}
\newcommand{\Ann}{\operatorname{Ann}}
\newcommand{\Rad}{\operatorname{Rad}}

\newcommand{\depth}{\operatorname{depth}}

\newcommand{\lo}{\longrightarrow}
\newcommand{\fm}{\frak{m}}
\newcommand{\fp}{\frak{p}}

\newcommand{\fa}{\frak{a}}

\begin{document}

\author[K. Divaani-Aazar and A. Hajikarimi]{Kamran Divaani-Aazar
and Alireza Hajikarimi}
\title[Generalized local cohomology modules and ...]
{Generalized local cohomology modules and homological Gorenstein
dimensions}

\address{K. Divaani-Aazar, Department of Mathematics, Az-Zahra
University, Vanak, Post Code 19834, Tehran, Iran-and-Institute for
Studies in Theoretical Physics and Mathematics, P.O. Box 19395-5746,
Tehran, Iran.} \email{kdivaani@ipm.ir}

\address{A. Hajikarimi, Science and Research Branch, Islamic Azad
University, Tehran, Iran.} \email{alihajikarimi@yahoo.com}

\subjclass[2000]{13D45, 13D05.}

\keywords{Artinian modules, attached prime ideals, cohomological
dimension, generalized local cohomology modules,
Gorenstein injective dimension, Gorenstein projective dimension.\\
The first author was supported by a grant from IPM (No. 86130114).}

\begin{abstract} Let $\fa$ be an ideal of a commutative
Noetherian ring $R$ and $M$ and $N$ two finitely generated $R$-modules.
Let $\cd_{\fa}(M,N)$ denote the supremum of the $i$'s such that
$H^i_{\fa}(M,N)\neq 0$. First, by using the theory of Gorenstein
homological dimensions, we obtain several upper bounds for
$\cd_{\fa}(M,N)$. Next, over a Cohen-Macaulay local ring $(R,\fm)$,
we show that $$\cd_{\fm}(M,N)=\dim R-\grade(\Ann_RN,M),$$ provided
that either projective dimension of $M$ or injective dimension of
$N$ is finite. Finally, over  such rings, we establish an analogue
of the Hartshorne-Lichtenbaum Vanishing Theorem in the context of
generalized local cohomology modules.
\end{abstract}

\maketitle

\section{Introduction}

Let $R$ be a commutative Noetherian ring with identity. The notion
of generalized local cohomology was introduced by Herzog in his
Habilitationsschrift \cite{He}. Let $\fa$ be an ideal of $R$
and $M$ and $N$ two $R$-modules. The $i$th generalized local cohomology
module of $M$ and $N$ with respect to $\fa$ is defined by
$H^i_{\fa}(M,N):=\underset{n}{\varinjlim}\Ext^{i}_{R}(M/\fa^{n}M,
N)$.

Henceforth, we assume that $M$ and $N$ are finitely generated. In
\cite[Proposition 5.5]{B}, it is shown that the least integer $i$
such that $H^i_{\fa}(M,N)\neq 0$ is equal to $\grade(\Ann_R(M/\fa
M),N)$. We denote the supremum of $i$'s such that
$H^i_{\fa}(M,N)\neq 0$ by $\cd_{\fa}(M,N)$ and we abbreviate
$\cd_{\fa}(R,N)$ by $\cd_{\fa}(N)$. In Section 2, we explore
interrelations between generalized local
cohomology modules and homological Gorenstein dimensions. This
direction of research was motivated by Sazeedeh's work
\cite{Sa1},  which implies that local cohomology modules can be
computed by Gorenstein injective resolutions. Here, we will apply
the theory of Gorenstein homological dimensions to establish the
following upper bounds for $\cd_{\fa}(M,N)$:
\begin{enumerate}
\item[i)] If $\Gpd_NM$ is finite, then $\cd_{\fa}(M,N)\leq
\Gpd_NM+\cd_{\fa}(M\otimes_RN)$.
\item[ii)] If $\pd_NM$ or $\id N$ is finite, then $\cd_{\fa}(M,N)
\leq \min\{\dim R,\Gpd_NM+\cd_{\fa}(M\otimes_RN)\}$.
\item[iii)] If $\Gpd M$ and $\pd N$ are finite,  then $\cd_{\fa}(M,N)
\leq \min\{\dim R,\Gpd_NM+\cd_{\fa}(M\otimes_RN)\}$.
\item[iv)] If $\id M$ and $\Gid N$ are finite, then $\cd_{\fa}(M,N)
\leq \min\{\dim R,\Gpd_NM+\cd_{\fa}(M\otimes_RN)\}$.
\item[v)] If either $\pd M$ or $\id M$ is finite, then $\cd_{\fa}
(M,N)\leq \min\{\Gid N,\Gpd_NM+\cd_{\fa}(M\otimes_RN)\}$.
\end{enumerate}
Here $\Gid$ and $\Gpd$ stand for Gorenstein injective and projective
dimensions, respectively. Also $\pd_NM:=\sup\{\pd_{R_{\fp}}M_{\fp}:
\fp\in \Supp_RM\cap \Supp_RN \}$ and $\Gpd_NM:=\sup\{i\in
\mathbb{N}_0:\Ext_R^i(M,N)\neq 0\}$, with the usual convention that
the supremum of the empty set of integers is interpreted as
$-\infty$. As an application of these bounds, one can improves
\cite[Theorem 3.2]{HZ} and Theorem 3.1 and Lemma 5.4 in \cite{CH},
see Corollary 2.6 below.

Let $(R,\fm,k)$ be a local ring. Grothendieck's non-Vanishing
Theorem implies that $\cd_{\fm}(N)=\dim N$. But, not much is known
about $\cd_{\fm}(M,N)$. In \cite{HZ}, the class of finitely
generated $R$-modules $L$ for which $\cd_{\fm}(L,N)=\depth N$ is
investigated. In this paper, we intend to compute $\cd_{\fm}(M,N)$,
when either projective dimension of $M$ or injective dimension of
$N$ is finite. Note that we have to impose these assumptions on $M$
or $N$, because otherwise $\cd_{\fm}(M,N)$ might be infinite. To
realize an easy example, suppose that $M$ has infinite projective
dimension. Then $H^i_{\fm}(M,k)\cong \Ext_R^i(M,k)$ is non-zero for
infinitely many $i$. Now, let $R$ be  Cohen-Macaulay and suppose
that either projective dimension of $M$ or injective dimension of $N$
is finite. In Section 3, we show that $$\cd_{\fm}(M,N)=\dim
R-\grade(\Ann_RN,M).$$  Next, we show that $H^d_{\fa}(M,N), d:=\dim R,$
can be described as a certain quotient of  $H^d_{\fm}(M,N)$, see Theorem
3.8 below. Then as an immediate application, we will compute the set of
attached prime ideals of the Artinian $R$-module $H^d_{\fa}(M,N)$.
More precisely, we show that
$$\Att_R(H^d_{\fa}(M,N))=\{\fp\cap R:\fp\in
\Supp_{\hat{R}}\hat{N}\cap \Ass_{\hat{R}}\hat{M} \ \ and \ \  \dim
(\hat{R}/\fa \hat{R}+\fp)=0\}.$$ This enables us to establish an analogue
of the Hartshorne-Lichtenbaum Vanishing Theorem for generalized local
cohomology modules to the effect that the following are
equivalent.
\begin{enumerate}
\item[i)] $H^d_{\fa}(M,N)=0$.
\item[ii)] $H^d_{\fm}(M,N)=\displaystyle{\sum_{n\in
\mathbb{N}}}<\fm>(0:_{H^d_{\fm}(M,N)}\fa^n)$.
\item[iii)] For any integer $l\in \mathbb{N}$, there exists an
$n=n(l)\in \mathbb{N}$ such that
$$0:_{H^d_{\fm}(M,N)}\fa^l\subseteq
<\fm>(0:_{H^d_{\fm}(M,N)}\fa^n).$$
\item[iv)] $\dim \hat{R}/\fa \hat{R}+\fp>0$ for all
$\fp\in \Supp_{ \hat{R}}\hat{N}\cap \Ass_{\hat{R}}\hat{M}$.
\end{enumerate}
(Here for an Artinian $R$-module $A$, we use $<\fm>A$ for denoting
the submodule $\underset{i\in \mathbb{N}}\cap \fm^i A$.)

\section{Some upper bounds for $\cd_{\fa}(M,N)$}

We start this section by introducing the notions of relative
projective dimension and of relative Gorenstein projective dimension
for a pair of finitely generated $R$-modules. But, first let recall
some definitions from the theory of Gorenstein homological
dimensions. The theory of Gorenstein homological dimensions was
initiated by Enochs and Jenda in 1995, see \cite{EJ}. An
$R$-module $N$ is said to be Gorenstein injective if there exists an
exact sequence
$$I^{\bullet}: \cdots \lo I_1\lo I_0\lo I^0\lo I^1\lo \cdots$$
of injective $R$-modules such that $N\cong \ker(I^0\lo I^1)$ and
that $\Hom_R(I,I^{\bullet})$ is exact for all injective $R$-modules
$I$. Also, an $R$-module $N$ is said to be Gorenstein projective if
there exists an exact sequence
$$P_{\bullet}: \cdots \lo P_1\lo P_0\lo P^0\lo P^1\lo \cdots$$
of projective $R$-modules such that $N\cong \ker(P^0\lo P^1)$ and
that $\Hom_R(P_{\bullet},P)$ is exact for all projective $R$-modules
$P$. Then for an $R$-module $N$, the notion $\Gid N$, Gorenstein
injective dimension of  $N$, is defined as the infimum of the length
of right resolutions of $N$ which are consisting of Gorenstein
injective modules. Similarly, $\Gpd N$, Gorenstein projective
dimension of $N$, is defined as the infimum of the length of left
resolutions of $N$ which are consisting of Gorenstein projective
modules. It is known that if $\id N$ (, respectively $\pd N$) is
finite, then $\Gid N=\id N$ (, respectively $\Gpd N=\pd N$).

By \cite[Corollary 2.21]{Ho}, if $\Gpd M<\infty$, then $\Gpd
M=\sup\{i\in \mathbb{N}_0:\Ext_R^i(M,R)\neq 0\}$. In view of this,
our second definition below might seem to be rational.

\begin{definition} Let $M$ and $N$ be two finitely generated
$R$-modules. We define {\it projective dimension of $M$ relative to
$N$} by $$\pd_NM:=\sup\{\pd_{R_{\fp}}M_{\fp}:\fp\in \Supp_RM\cap
\Supp_RN \}.$$ Also, we define {\it Gorenstein projective dimension
of $M$ relative to $N$} by
$$\Gpd_NM:=\sup\{i\in \mathbb{N}_0:\Ext_R^i(M,N)\neq 0 \}.$$
\end{definition}

Note that the above mentioned result of Holm indicates that if
$\Gpd M<\infty$, then $\Gpd_RM=\Gpd M$. For sake of completeness, we
collect all needed properties of these newly defined notions in a lemma.

\begin{lemma} Let $M$ and $N$ be two non-zero finitely generated
$R$-modules.
\begin{enumerate}
\item[i)] $\grade(\Ann _R M,N)\leq \Gpd_NM$.
\item[ii)] $\pd_NM\leq \pd M$ and if $R$ is local, then $\pd_NM=
\pd M$.
\item[iii)] $\Gpd_NM\leq \min \{\pd M, \id N \}$. In particular,
if either $\pd M$ or $\id N$ is finite, then $\Gpd_NM$ is finite.
\item[iv)] If $\pd_NM$ is finite, then
$\Gpd_NM=\pd_NM$. In particular, if $R$ is local and $\pd M$ is
finite, then $\Gpd_NM=\pd M$.
\item[v)] If $\Gpd M$ and $\pd N$ are finite, then $\Gpd_NM\leq
\Gpd M$.
\item[vi)] If $\id M$ and $\Gid N$ are finite, then $\Gpd_NM\leq
\Gid N$.
\item[vii)] If $\pd_NM$ is finite, then $\bigcup_{i\in \mathbb{N}_0}
\Supp_R (\Ext_R^i(M,N))=\Supp_R(M\otimes_RN)$, and so  $\dim
(M\otimes_RN)=\max\{\dim \Ext_R^i(M,N):i\in \mathbb{N}_0 \}$.
\item[viii)] If $\pd_NM$ is finite, then for any ideal $\fa$ of $R$,
$\cd_{\fa}(M\otimes_RN)=\max\{\cd_{\fa}(\Ext_R^i(M,N)):i\in
\mathbb{N}_0 \}$.
\item[ix)] If $R$ is local and $\id
N<\infty$, then $\Gpd_NM=\depth R-\depth M$. In addition, if $M$ is
maximal Cohen-Macaulay, then $\Gpd_NM=0$.
\item[x)] If $\id N<\infty$, then $\Gpd_NM=\sup\{\Ht\fp-
\depth M_{\fp}:\fp\in \Supp_RM\cap\Supp_RN\}$.
\item[xi)] $\id_{R_{\fp}}N_{\fp}<\infty$ for all prime ideals
$\fp$ of $R$ if and only if $\Gpd_NL<\infty$ for all finitely
generated $R$-modules $L$.
\end{enumerate}
\end{lemma}

{\bf Proof.} i), ii) and iii) follow immediately by the definitions.

iv) Let $i>\pd_NM$. Then for any $\fp\in \Supp_RM\cap \Supp_RN$, it
turns out that  $$\Ext_{R}^{i}(M,N)_{\fp}\cong
\Ext_{R_{\fp}}^{i}(M_{\fp},N_{\fp})=0,$$ as
$i>\pd_{R_{\fp}}M_{\fp}$. Hence $\Gpd_{N}M\leq \pd_NM$.

Over a local ring $T$, for any two non-zero finitely generated
$T$-modules $M$ and $L$, \cite[page 154, Lemma 1]{Mat} implies that if $\pd
M<\infty$, then $\Gpd_LM=\pd M$ . Let $\fp\in \Supp_RM\cap \Supp_RN$
and $h:=\pd_{R_{\fp}}M_{\fp}$. Then
$$\Ext_R^h(M,N)_{\fp}\cong
\Ext_{R_{\fp}}^{h}(M_{\fp},N_{\fp})\neq 0,$$ and so
$\pd_{R_{\fp}}M_{\fp}\leq \Gpd_{N}M$. Therefore $\Gpd_{N}M=\pd_NM$.

v) and vi) follow by \cite[Theorem 2.20]{Ho} and \cite[Theorem
2.22]{Ho}, respectively.

vii) It is obvious that each module $\Ext_R^i(M,N)$ is supported in
$\Supp_RM\cap \Supp_RN=\Supp_R(M\otimes_RN)$. Now, let $\fp\in
\Supp_RM\cap \Supp_RN$ and set $h:=\pd_{R_{\fp}}M_{\fp}$. We
observed earlier in the proof of part iv) that $\Ext_{R_{\fp}}^h
(M_{\fp},N_{\fp})\neq 0$. Hence $\fp\in\Supp_R(\Ext_{R}^{h}(M,N))$,
and the conclusion follows.

viii) For any two finitely generated $R$-modules $X$ and $Y$ with
$\Supp_RX\subseteq \Supp_RY,$ \cite[Theorem 2.2]{DNT} implies that
$\cd_{\fa}(X)\leq \cd_{\fa}(Y)$. Let $Z$ be a finitely generated
$R$-module. Since $Z$ and $\oplus_{\fp\in \Ass_RZ} R/\fp$ have the
same support, it follows that
$$\cd_{\fa}(Z)=\sup\{\cd_{\fa}(R/\fp): \fp\in \Supp_RZ\}.$$ Now,
the claim becomes clear in light of part vii).

ix) By \cite[2.6]{I}, it follows that $\sup\{i\in
\mathbb{N}_0:\Ext_R^i(M,N)\neq 0 \}=\depth R-\depth M.$

x) Let $\fp\in \Supp_RN$. Since $\id_RN<\infty$, it turns out that
$\id_{R_{\fp}}N_{\fp}< \infty$ , and so $R_{\fp}$ is Cohen-Macaulay.
Hence, it follows by ix) that
$$\begin{array}{ll} \Gpd_NM&=
\sup\{\Gpd_{N_{\fp}}M_{\fp}:\fp\in \Supp_RM\cap \Supp_RN\}\\
&=\sup\{\Ht \fp-\depth_{R_{\fp}}M_{\fp}:\fp\in \Supp_RM\cap
\Supp_RN\}.
\end{array}$$

xi) is the main result of \cite{T}. $\Box$

\begin{example} There exist a non-Gorenstein Artinian local ring
$(R,\fm)$ and a finitely generated $R$-module $M$, which is not
Gorenstein projective such that $\Hom_R(M,R)\neq 0$ and
$\Ext_R^i(M,R)=0$ for all $i>0$. So, although neither projective
dimension of $M$ nor injective dimension of $R$ is finite, one has
$\Gpd_RM=0<\infty$. For a concrete realization of this example, we
refer the reader to \cite[Theorem 1.7]{JS}.
\end{example}

We need the following lemma in the proof of Theorem 2.5 and Corollary
2.6 below. It generalizes \cite[Lemma 3.1]{HZ}.

\begin{lemma} Assume that $M$ and $N$ are non-zero finitely
generated $R$-modules such that $\Gpd_NM<\infty$. Then obviously $\dim
\Ext_R^i(M,N)\leq \dim (M\otimes_RN)+\Gpd_NM-i$ for all $i$.
Moreover, if the assumption $\Gpd_NM<\infty$ replaced by any of the
following stronger assumptions, then $\dim \Ext_R^i(M,N)\leq \min
\{\dim R,\dim (M\otimes_RN)+\Gpd_NM\}-i$ for all $i$.
\begin{enumerate}
\item[i)] $\Gpd M$ and $\pd N$ are finite.
\item[ii)] $\id M$ and $\Gid N$ are finite.
\item[iii)] $\id N$ is finite.
\item[iv)] $\pd_NM$ is finite.
\end{enumerate}
\end{lemma}

{\bf Proof.} Let $i\in \mathbb{N}_0$. If $\Ext_R^i(M,N)=0$, then
there is nothing to prove. So, we may assume that $i\leq \Gpd_NM$.
On the other hand, $\dim \Ext_R^i(M,N)\leq \dim (M\otimes_RN)$,
since $\Ext_R^i(M,N)$ is supported in $\Supp_R(M\otimes_RN)$. Hence
$\dim \Ext_R^i(M,N)+i\leq \dim (M\otimes_RN)+\Gpd_NM$.

For a while, we assume that $R$ is local. By \cite[3.7]{AB}, if
$\Gpd M<\infty$, then $\Gpd M=\depth R-\depth M$. By \cite{KTY} (,
respectively \cite[Theorem 18.9]{Mat}), if $\Gid N$ (, respectively
$\id N$) is finite, then $\Gid N$ (, respectively $\id N$) is equal
to $\depth R$. Also, by the Auslander-Buchsbaum formula if $\pd
M<\infty$, then $$\Gpd_NM=\pd M=\depth R-\depth M.$$ Hence in view
of Lemma 2.2, each of the conditions   i), ii), iii) and iv) implies
that $\Gpd_NM\leq \depth R$. But for any $\fp\in \Supp_RM\cap
\Supp_RN$, all of the conditions  i), ii), iii) and iv) are
preserved under localization at $\fp$. Thus, for any such prime
ideal $\fp$, each of the conditions i), ii), iii) and iv) implies
that $\Gpd_{N_{\fp}}M_{\fp}\leq \Ht \fp$.

Assume one of the conditions i), ii), iii) and iv) holds and let
$\fp\in \Supp_R(\Ext_R^i(M,N))$. Clearly, we can suppose that $\dim
R$ is finite. Then we have $i\leq \Gpd_{N_{\fp}}M_{\fp}\leq \Ht
\fp$, and so $\dim R/\fp\leq \dim R-\Ht \fp\leq \dim R-i$.  Thus
$$\dim \Ext_R^i(M,N)=\sup\{\dim R/\fp:\fp\in
\Supp_R(\Ext_R^i(M,N))\}\leq \dim R-i.$$ Now, in view of the first
assertion, the proof is complete. $\Box$

Let $\fa$ be an ideal of $R$ and $M$ a finitely generated $R$-module
of finite dimension $d$. It is known that $H^d_{\fa}(M)$ is
Artinian. Also, by \cite[Corollary 2.5]{Mar}, we know that
$\Supp_RH^{d-1}_{\fa}(M)$ is finite. The second assertion of the
following result might be considered as a generalization of these
facts to generalized local cohomology modules.

\begin{theorem} Let $\fa$ be an ideal of $R$ and $M$ and $N$ two
finitely generated $R$-modules such that $p:=\Gpd_NM<\infty$ and set
$c:=\Gpd_NM+\cd_{\fa}(M\otimes_RN)$. Then $H^c_{\fa}(M,N)\cong
H^{c-p}_{\fa}(\Ext^p_R(M,N))$ and $H^i_{\fa}(M,N)=0$ for all $i>c$.
Moreover, if $\dim (M\otimes_RN)<\infty$ and $d:=\dim
(M\otimes_RN)+\Gpd_NM$, then $H^d_{\fa}(M,N)$ is Artinian and
$\Supp_R(H^{d-1}_{\fa}(M,N))$ is finite.
\end{theorem}

{\bf Proof.} First of all, we claim that $H_{\fa}^i(\Hom_R(M,E))=0$ for any
injective $R$-module $E$ and all $i\geq 1$. Since any injective $R$-module decomposes into
a direct sum of indecomposable injective $R$-modules, we may and do assume that
$E=E_R(R/\fp)$, for some prime ideal $\fp$ of $R$. (Note that the functor $H_{\fa}^i(\cdot)$ commutes with
direct sums, and as $M$ is finitely generated the functor $\Hom_R(M,\cdot)$ also commutes with direct sums.)
Since $$\Hom_R(M,E)\cong \Hom_R(M,\Hom_{R_{\fp}}(R_{\fp},E))\cong
\Hom_{R_{\fp}}(M_{\fp},E)\cong \Hom_{R_{\fp}}(M_{\fp},E_{R_{\fp}}(R_{\fp}/\fp R_{\fp})),$$ we
deduce that $\Hom_R(M,E)$ is an Artinian $R_{\fp}$-module. Hence
$H_{\fa}^i(\Hom_R(M,E))=H_{\fa R_{\fp}}^i(\Hom_R(M,E))=0$, as claimed.
Let $F(\cdot):=\Gamma_{\fa}(\cdot)$ and $G(\cdot):=\Hom_R(M,\cdot)$.  Since
$(GF)(\cdot)=\Hom_R(M,\Gamma_{\fa}(\cdot))$, by
\cite [Theorem 11.38]{R}, one has the following Grothendieck's spectral sequence
$$E^{i,j}_2:=H^{i}_{\fa} (\Ext_R^j(M,N))\underset{i}\Longrightarrow
H^{i+j}_{\fa}(M,N).$$ Hence for each $n\in\mathbb{N}_{0}$, there
exists a chain $$0=H^{-1}\subseteq H^{0}\subseteq\cdots \subseteq
H^{n}:=H^{n}_{\fa} (M,N) \  \  (\ast) $$ of submodules of
$H^{n}_{\fa} (M,N)$ such that $H^{i}/H^{i-1}\cong
E^{i,n-i}_{\infty}$ for all $i=0,1,\cdots,n$.

If $n>\Gpd_{N}(M)+\cd_{\fa}(M\otimes_RN)$, then either
$i>\cd_{\fa}(M\otimes_RN)$ or $n-i>\Gpd_{N}M$. In each case,
$E_{\infty}^{i,n-i}=0$, as $E^{i,n-i}_{\infty}$ is a subquotient of
$E^{i,n-i}_{2}$. (Note that  for each $j$,
$\cd_{\fa}(\Ext_R^j(M,N))\leq \cd_{\fa}(M\otimes_RN)$.) Therefore,
from the chain $(\ast)$, it follows that $H^{n}_{\fa}(M,N)=0$  for
all $n>\Gpd_{N}(M)+\cd_{\fa}(M\otimes_RN)$.

Next, in the chain $(*)$, let $n=c$. Since $E_2^{i,j}=0$, whenever
$i>c-p$ or $j>p$, it turns out that $H^{c}_{\fa}(M,N)\cong
E_{\infty}^{c-p,p}$. For each $r\geq 2$, consider the sequence
$$E_r^{c-p-r,p+r-1} \overset{d_r^{c-p-r,p+r-1}} \lo E_r^{c-p,p}
\overset{d_r^{c-p,p}} \lo E_r^{c-p+r,p-r+1}.$$ Since
$E_r^{c-p-r,p+r-1}$ and $E_r^{c-p+r,p-r+1}$ are respectively
subquotients of $E_2^{c-p-r,p+r-1}$ and $E_2^{c-p+r,p-r+1}$, it
follows that $E_r^{c-p-r,p+r-1}=E_r^{c-p+r,p-r+1}=0$, and so
$$E_{r+1}^{c-p,p}=\ker d_r^{c-p,p}/\im d_r^{c-p-r,p+r-1}\cong
E_r^{c-p,p}.$$ Hence $H^{c}_{\fa}(M,N)\cong E_{\infty}^{c-p,p}\cong
\dots \cong E_2^{c-p,p}$.

Now, assume that $\dim(M\otimes_RN)<\infty$ and set
$d:=\dim(M\otimes_RN)+\Gpd_NM$. Clearly, $c\leq d$ and in view of
the first assertion of the theorem, we may and do assume that $d=c$.
As $H^c_{\fa}(M,N)\cong H^{c-p}_{\fa}(\Ext^p_R(M,N))$ and $\dim
\Ext^p_R(M,N)\leq c-p$, it follows that $H^c_{\fa}(M,N)$ is
Artinian.

Let $i,j$ be two non-negative integers such that $i+j=d-1$. Lemma
2.4 yields that $\dim \Ext^j_R(M,N)\leq d-j$, and so $i=d-j-1\geq
\dim \Ext^j_R(M,N)-1$. So, by \cite[Corollary 2.5]{Mar},
$\Supp_R(E_{\infty}^{i,j})$ is finite. Thus, from the chain $(*)$,
we deduce that $\Supp_R(H^{d-1}_{\fa}(M,N))$ is finite. $\Box$

In the above argument, we used the first assertion of Lemma 2.4.
Slightly modifying it in the light of the second assertion of Lemma
2.4, concludes the following corollary. It is worth to mention that
this corollary improves Theorem 3.1 and Lemma 5.4 in \cite{CH}.
It also generalizes \cite[Theorem 3.2]{HZ}.

\begin{corollary} Let $\fa$ be an ideal of $R$ and $M$ and $N$ two
finitely generated $R$-modules. Assume that one of the following conditions
is satisfied:
\begin{enumerate}
\item[i)] $\Gpd M$ and $\pd N$ are finite.
\item[ii)] $\id M$ and $\Gid N$ are finite.
\item[iii)] $\id N$ is finite.
\item[iv)] $\pd_NM$ is finite.
\end{enumerate}
Then $H^i_{\fa}(M,N)=0$ for all $i>\min\{\dim R,
\Gpd_NM+\cd_{\fa}(M\otimes_RN)\}$. Moreover, if $\dim
(M\otimes_RN)<\infty$ and $d:=\min\{\dim R,\dim
(M\otimes_RN)+\Gpd_NM\}$, then $H^d_{\fa}(M,N)$ is Artinian and
$\Supp_R(H^{d-1}_{\fa}(M,N))$ is finite.
\end{corollary}

Among other things, Theorem 2.5 says that $\cd_{\fa}(M,N)\leq
\Gpd_NM+\cd_{\fa}(M\otimes_RN)$. As the following corollary
indicates, in some cases the equality holds.

\begin{corollary} Let $\fa$ be an ideal of $R$ and $M$ and $N$ two
finitely generated $R$-modules.
\begin{enumerate}
\item[i)]If $p:=\Gpd_N M=\grade(\Ann_RM,N)$, then $H^i_{\fa}
(M,N)=H^{i-p}_{\fa}(\Ext^p_R(M,N))$ for all i, and consequently
$\cd_{\fa}(M,N)=\Gpd_NM+\cd_{\fa}(\Ext^p_R(M,N))$.
\item[ii)] If $p:=\pd_N M=\grade(\Ann_RM,N)$, then $\cd_{\fa}(M,N)
=\Gpd_NM+\cd_{\fa}(M\otimes_R N)$.
\item[iii)] If $R$ is local, $\id N<\infty$ and $M$  maximal
Cohen Macaulay and faithful, then $\cd_{\fa}(M,N)=\cd_{\fa}(N)$.
\end{enumerate}
\end{corollary}

{\bf Proof.} i) The spectral sequence
$H^{i}_{\fa}(\Ext_{R}^{j}(M,N))\underset{i}\Longrightarrow H^{i+j}
_{\fa} (M,N)$ collapses  at $j=p$, and so $H^n _\fa (M,N)\cong
H^{n-p}_{\fa}(\Ext_R^p(M,N))$ for all $n$. This shows that
$\cd_{\fa}(M,N)=p+\cd_{\fa}(\Ext^p_R(M,N))$.

ii) Lemma 2.2 viii) yields that
$\cd_{\fa}(\Ext_R^p(M,N))=\cd_{\fa}(M\otimes_RN)$. Hence ii) follows
by i).

iii) Lemma 2.2 ix) yields that $\Gpd_NM=0$, and so by the proof of
i), we have the isomorphisms $H^n_{\fa}
(M,N)=H^n_{\fa}(\Hom_R(M,N))$ for all $n$. Since by \cite[Exercise
1.2.27]{BH}, $\Ass_R(\Hom_R(M,N))= \Supp_RM\cap \Ass_RN$, it follows
that $\Supp_R(\Hom_R(M,N))=\Supp_RN.$ Now, \cite[Theorem 2.2]{DNT}
yield that
$$\cd_{\fa}(M,N)=\cd_{\fa}(\Hom_R(M,N))=\cd_{\fa}(N). \  \ \Box$$

For two finitely generated $R$-modules $M$ and $N$, Theorem 2.5
implies that if $\Gpd_NM<\infty$, then $\cd_{\fa}(M,N)<\infty$ for
all ideals $\fa$ of $R$. The second part of the next corollary
indicates that the converse is also true.

\begin{corollary} Let $M$ and $N$ be two finitely generated $R$-modules.
\begin{enumerate}
\item[i)] Let $\fa$ be an ideal of $R$ such that $\Supp_RM\cap
\Supp_RN\subseteq \V(\fa)$. Then $H^i_{\fa}(M,N)\cong \Ext_R^i(M,N)$
for all $i\in \mathbb{N}_0$.
\item[ii)] $\Gpd_NM<\infty$ if and only if $H^i_{\fa}(M,N)=0$ for all
ideals $\fa$ of $R$ and for all $i\gg 0$.
\end{enumerate}
\end{corollary}

{\bf Proof.} i) Since $$\Supp_R(\Ext_R^j(M,N))\subseteq \Supp_RM\cap
\Supp_RN\subseteq \V(\fa),$$ it follows that $\Ext_R^j(M,N)$ is
$\fa$-torsion for all $i$. Hence, the spectral sequence
$$H^{i}_{\fa}(\Ext_{R}^{j}(M,N))\underset{i}\Longrightarrow
H^{i+j} _{\fa} (M,N)$$ collapses  at $i=0$, and so $$H^n_{\fa}
(M,N)\cong H^0_{\fa}(\Ext_R^n(M,N))=\Ext_R^n(M,N)$$ for all $n$.

ii) The ``only if'' part follows by Theorem 2.5. For the converse,
let $\fa:=\Ann_R(M\otimes_RN)$. Then by i), we have
$\Ext_R^i(M,N)\cong H_{\fa}^i(M,N)=0$ for all $i\gg 0$. Hence
$\Gpd_NM<\infty$. $\Box$

Sazeedeh \cite{Sa1} has proved that local cohomology modules of
an $R$-module $N$ can be computed by Gorenstein injective
resolutions of $N$. Also, he \cite{Sa2} proved that if $\dim R$ and
$\pd M$ are finite, then $H^i_{\fa}(M,N)$ can be computed by using
Gorenstein injective resolutions of $N$. Next, by employing his
technique that was used in the proof of \cite[Theorem 3.1]{Sa1}, we
improve his later-mentioned result.

\begin{lemma} Let $\fa$ be an ideal of $R$ and $M$ a finitely
generated $R$-module such that either $\pd M$ or $\id M$ is finite.
Then $H^{i}_{\fa}(M,N)=0$ for any Gorenstein injective $R$-module
$N$ and all $i>0$. Hence for any $R$-module $N$, the generalized
local cohomology modules $H^{i}_{\fa}(M,N)$ can be computed by
Gorenstein injective resolutions of $N$.
\end{lemma}

{\bf Proof.}  Let $N$ be a Gorenstein injective $R$-module. Let
$c:=\min\{\pd M,\id M\}+\ara(\fa)$ and $E^{i,j}_2:=H^{i}_{\fa}(
\Ext_R^j(M,N))$ for all $i,j\geq 0$.  If $\id M<\infty$, then by
Lemma 2.2 vi), $\Gpd_NM=0$. Therefore, the spectral sequence
$E^{i,j}_2\underset{i}\Longrightarrow H^{i+j} _{\fa}(M,N)$ implies
that $H^{i}_{\fa}(M,N)=0$ for all $i>c$ if either $\pd M<\infty$ or
$\id M<\infty$.  Since $N$ is Gorenstein injective, there exists an
exact sequence $$I^{\bullet}: \cdots \lo I_1\lo I_0\lo I^0\lo I^1\lo
\cdots$$ of injective $R$-modules such that $N\cong \ker(I^0\lo
I^1)$ and that $\Hom_R(I,I^{\bullet})$ is exact for all injective
$R$-modules $I$. For each $i> 0$, let $N_i:= \ker(I_i\lo I_{i-1})$,
$N_0:=\ker(I_0\lo I^0)$ and $N_{-1}:=N$. For each $i\geq 0$, from
the exact sequence $$0\lo N_i\lo I_i\lo N_{i-1}\lo 0,$$ we deduce
the following long exact sequences of generalized local cohomology
modules
$$\cdots \lo H_{\fa}^{j}(M,I_i)\lo H_{\fa}^{j}(M,N_{i-1})\lo
H_{\fa}^{j+1}(M,N_i)\lo H_{\fa}^{j+1}(M,I_i)\lo \cdots.$$ Thus, we
conclude the isomorphisms $H^j_{\fa}(M,N_{i-1})\cong
H^{j+1}_{\fa}(M,N_i)$ for all $i\geq 0$ and all $j\geq 1$. Hence for
$i>0$, we have
$$H^i_{\fa}(M,N)\cong H^{i+1}_{\fa}(M,N_0)\cong \dots
\cong H^{i+c}_{\fa}(M,N_{c-1})=0.$$ (Note that $N_{c-1}$ is
Gorenstein injective.) Thus any Gorenstein injective $R$-module is
$H^0_{\fa}(M,\cdot)$-acyclic. This finishes the proof. Recall that
if $T$ is a left exact additive functor from the category of
$R$-modules and $R$-homomorphisms to itself, then for any $R$-module
$N$, the right derived functors $R^iT$ of $T$ at $N$ can be computed
by using right resolutions of $N$ which are consisting of T-acyclic
modules. $\Box$

\begin{corollary} Let $\fa$ be an ideal of $R$ and $M,N$  finitely
generated $R$-modules such that either $\pd M$ or $\id M$ is finite.
Then $\cd_{\fa}(M,N)\leq \min\{\Gid
N,\Gpd_NM+\cd_{\fa}(M\otimes_RN)\}$.
\end{corollary}

\begin{example}  When $\Gpd M$ is finite, it is rather natural to ask
whether $H_{\fa}^i(M,N)$ can be computed by Gorenstein injective
resolutions of $N$.
This would not be the case. To see this, let $(R,\fm,k)$ be a
Gorenstein local ring which is not regular. Then $\Gpd k$ and $\Gid
k$ are both finite, while $H_{\fm}^i(k,k)\cong \Ext_R^i(k,k)$ is
non-zero for infinitely many $i$.
\end{example}

\section{The local case}

In this section, for a pair $(M,N)$ of finitely generated modules
over a Cohen-Macaulay local ring $(R,\fm)$, we compute
$\cd_{\fm}(M,N)$ provided that either $\pd M$ or $\id N$ is finite.
We will end this section by proving an analogue of the
Hartshorne-Lichtenbaum Vanishing Theorem in the context of generalized local
cohomology modules. The main ingredients in our proofs are Suzuki's
and the Herzog-Zamani Duality Theorems for generalized local
cohomology modules. We start this section with the following useful result.

\begin{lemma} Let $\fa$ be an ideal of $R$ an $x$ an element of $R$.
Let $M$ be a finitely generated $R$-module and $N$ an arbitrary
$R$-module. There is a natural long exact sequence
$$\cdots\lo H_{\fa+(x)}^i(M,N)\lo H_{\fa}^i(M,N)\lo H_{\fa R_x}^
i(M_x,N_x)\lo H_{\fa+(x)}^{i+1}(M,N)\lo \cdots .$$
\end{lemma}

{\bf Proof.} Let $$I^{\bullet}: 0\lo I^0\overset{d^0}\lo I^1\lo
\cdots \lo I^i\overset{d^i}\lo I^{i+1}\lo \cdots $$ be an injective
resolution of $N$. Since $M$ is finitely generated, \cite[Lemma 2.1
i)]{DST} implies that $H_{\fa}^i(M,N)\cong
H^i(\Hom_R(M,\Gamma_{\fa}(I^{\bullet})))$. Similarly, since the
$R_x$-module $M_x$ is finitely generated and $I^{\bullet}_x$, the
localization of $I^{\bullet}$ at $x$, provides an injective
resolution for the $R_x$-module $N_x$, it turns out that $H_{\fa
R_x}^i(M_x,N_x)\cong H^i(\Hom_{R_x}(M_x,\Gamma_{\fa
R_x}(I^{\bullet}_x)))$.

Let $I$ be an injective $R$-module. \cite[Lemma 8.1.1]{BS} yields
the following split exact sequence
$$0\lo \Gamma_{\fa+(x)}(I)\overset{i}\lo \Gamma_{\fa}(I)\overset{f}
\lo \Gamma_{\fa}(I_x)\lo 0,$$ where the maps are the natural ones.
Because of the natural isomorphism
$$\Hom_{R_x}(M_x,\Gamma_{\fa R_x}(I_x))\cong
\Hom_R(M,\Gamma_{\fa}(I_x)),$$ we can deduce the following exact
sequence of complexes
$$0\lo \Hom_R(M,\Gamma_{\fa+(x)}(I^{\bullet}))\overset{i^{\bullet}}
\lo \Hom_R(M,\Gamma_{\fa}(I^{\bullet})) \overset{f^{\bullet}}\lo
\Hom_{R_x}(M_x,\Gamma_{\fa R_x}(I^{\bullet}_x))\lo 0.$$ Its long
exact sequence of cohomologies is precisely our desired long exact
sequence. $\Box$

Henceforth, we assume $(R,\fm)$ is a local ring and $M, N$ finitely
generated $R$-modules. We intend to compute $\cd_{\fm}(M,N)$. First,
we specialize Theorem 2.5 and Lemma 2.6 to the case $\fa=\fm$.

\begin{corollary} Let $(R,\fm)$ be a local ring and $M$ and $N$ two
finitely generated $R$-modules. Consider the following conditions:
\begin{enumerate}
\item[i)] $\Gpd_NM<\infty.$
\item[ii)] $\Gpd M$ and $\pd N$ are finite.
\item[iii)] $\id M$ and $\Gid N$ are finite.
\item[iv)] $\id N$ is finite.
\item[v)] $\pd_NM$ is finite.
\end{enumerate}
Put $d:=\dim (M\otimes_RN)+\Gpd_NM$ in case i) and $d:=\min\{\dim
R,\dim (M\otimes_RN)+\Gpd_NM\}$ in other cases. In each of the above
cases, $\cd_{\fm}(M,N)\leq d$. Moreover for any ideal $\fa$ of $R$,
in each of the above cases, $H_{\fa}^d(M,N)$ is a homomorphic image
of $H_{\fm}^d(M,N)$.
\end{corollary}

{\bf Proof.} By Grothendieck's  non-Vanishing Theorem
$\cd_{\fm}(M\otimes _RN)=\dim (M\otimes_RN)$. Hence, the first
assertion is immediate by Theorem 2.5 and Corollary 2.6.

Now, we prove the second assertion. We may choose  $x_1,x_2, \dots
,x_n \in R$ such that $\fm=\fa+(x_1,x_2,\dots ,x_n)$. Set $\fa_i:
=\fa+(x_1,\dots ,x_{i-1})$ for $i=1,\dots , n+1$. For each $1\leq
i\leq n$, by Lemma 3.1, we have the following long exact sequence of
generalized local cohomology modules
$$\cdots
\lo H_{\fa_{i+1}}^d(M,N)\lo H_{\fa_i}^d(M,N)\lo H_{\fa_i R_{x_i}}
^d(M_{x_i},N_{x_i})\lo \cdots .$$ By Theorem 2.5 or Corollary 2.6,
$H_{\fa_i}^d(M,N)$ is Artinian. Hence $H_{\fa_i}^d(M,N)$ is
supported at most at $\fm$, and so $$H_{\fa_i
R_{x_i}}^d(M_{x_i},N_{x_i})\cong H_{\fa_i}^d(M,N) _{x_i}=0.$$ Hence
the natural homomorphism $H_{\fa_{i+1}}^d(M,N)\lo H_{\fa_i}^d(M,N)$
is epic. Using this successively for $1\leq i\leq n$, yields that
$H_{\fa}^d(M,N)$ is a homomorphic image of $H_{\fm}^d(M,N)$. $\Box$

The key to the proof of Theorem 3.5 below is given in the following
lemma.

\begin{lemma} Let $(R,\fm,k)$ be a Cohen Macaulay local ring and
$M,N$  finitely generated $R$-modules such that $\pd M<\infty$. Let
$\omega_{\hat{R}}$ denote the canonical module of $\hat{R}$.
\begin{enumerate}
\item[i)] $\grade(\fa \hat{R},\hat{M}\otimes_{\hat{R}}
\omega_{\hat{R}})=\grade(\fa,M)$ for any ideal $\fa$ of $R$.
\item[ii)] $\Ass_R(\hat{M}\otimes_{\hat{R}}\omega_{\hat{R}})=
\Ass_RM$.
\item[iii)] $\Rad(\Ann_{\hat{R}}(\Hom_{\hat{R}}(\omega_{\hat{R}},
\hat{N})))=\Rad(\Ann_{\hat{R}}\hat{N})$. Moreover, if $\id
N<\infty$, then $\Ann_{\hat{R}}(\Hom_{\hat{R}}(\omega_{\hat{R}},
\hat{N}))=\Ann_{\hat{R}}\hat{N}$.
\end{enumerate}
\end{lemma}

{\bf Proof.} i) We first prove the claim for $\fa=\fm$. Let $d:=\dim
R$. Then, by Suzuki's Duality Theorem
$$\Ext_{\hat{R}}^i(k,\hat{M}\otimes_{\hat{R}}\omega_{\hat{R}})\cong
H_{\fm}^{d-i}(M,k)^{\vee} \cong \Ext_R^{d-i}(M,k)^{\vee}.$$ Hence,
it follows by the Auslander-Buchsbaum formula that
$$
\begin{array}{ll} \depth (\hat{M}\otimes_{\hat{R}}\omega_{\hat{R}})&
=\inf\{i:\Ext_R^{d-i}(M,k)^{\vee}\neq 0\}\\
&=\inf\{d-j:\Ext_R^j(M,k)\neq 0\}\\
&=d-\sup\{j:\Ext_R^j(M,k)\neq 0\}\\
&=\depth M.
\end{array}$$
In particular, if $R$ possesses a canonical module $\omega_R$, then
$\depth M=\depth (M\otimes_R\omega_R)$.

Now, we prove the claim for an arbitrary ideal $\fa$. Without loss
of generality, we may assume that $R$ is complete. It follows by
\cite[Theorem 3.3.5 b)]{BH}, that for each prime ideal $\fp$, the
$R_{\fp}$-module $(\omega_R)_{\fp}$ is the canonical module of
$R_{\fp}$. For any finitely generated $R$-module $L$, by
\cite[Proposition 1.2.10 a)]{BH}, we have
$\grade(\fa,L)=\inf\{\depth L_{\fp}:\fp\in \V(\fa)\}$. Thus
$$\grade(\fa,M)=\inf \{\depth
(M_{\fp}\otimes_{R_{\fp}}(\omega_R)_{\fp}):\fp\in
\V(\fa)\}=\grade(\fa,M\otimes_R \omega_R).$$

ii) Since the set $\Ass_R(\hat{M}\otimes_{\hat{R}}\omega_
{\hat{R}})$ (, respectively $\Ass_RM$) is precisely consisting of
the contractions of elements of $\Ass_{\hat{R}}(\hat{M}\otimes_
{\hat{R}}\omega_{\hat{R}})$ (, respectively $\Ass_{\hat{R}}
\hat{M}$) to $R$, we can assume that $R$ is complete. Remember that
by the first paragraph of the proof, $\depth M_{\fp}=\depth
(M\otimes_R\omega_R)_{\fp}$ for all prime ideal $\fp$ of $R$. But,
for a finitely generated $R$-module $L$, one can check easily that
$\fp\in \Ass_RL$ if and only if $\depth L_{\fp}=0$. This yields that
$\Ass_R(M\otimes_R\omega_{R})=\Ass_RM$.

iii) Since by \cite[Theorem 3.3.5 b)]{BH},
$\Supp_{\hat{R}}\omega_{\hat{R}}=\Spec \hat{R}$, it turns out by
\cite[Exercise 1.2.27]{BH} that
$$\Ass_{\hat{R}}(\Hom_{\hat{R}}(\omega_{\hat{R}},\hat{N}))=
\Supp_{\hat{R}}\omega_{\hat{R}}\cap
\Ass_{\hat{R}}\hat{N}=\Ass_{\hat{R}}\hat{N}.$$ Consequently,
$\Rad(\Ann_{\hat{R}}(\Hom_{\hat{R}}(\omega_{\hat{R}},
\hat{N})))=\Rad(\Ann_{\hat{R}}\hat{N})$, as required. If $\id
N<\infty$, then by \cite[Theorem 4.3 ii)]{Su}, there is a natural
isomorphism $\hat{N}\cong \Hom_{\hat{R}}(\omega_{\hat{R}},
\hat{N})\otimes_{\hat{R}}\omega_{\hat{R}}$, which clearly implies
that $\Ann_{\hat{R}}(\Hom_{\hat{R}}(\omega_{\hat{R}},
\hat{N}))=\Ann_{\hat{R}}\hat{N}$. $\Box$

The statement of the corollary below involves the notion of attached
prime ideals. For convenient of the reader, we review this notion
briefly in below. Let $A=A_1+\dots +A_n$ be a minimal secondary
representation of the Artinian $R$-module $A$. Then the ideals
$\fp_i:=\Rad(\Ann_RA_i)$'s are prime and they are independent of the
given minimal secondary representation. Each $\fp_i$ is said to be
an attached prime ideal of $A$ and the set $\{\fp_1,\dots , \fp_n\}$
is denoted by $\Att_RA$. It is easy to see that a prime ideal $\fp$
of $R$ is an attached prime ideal of $A$ if and only if there exists
a quotient $C$ of $A$ such that $\fp=\Ann_RC$. In particular, this
implies that $A$ is zero if and only if $\Att_RA$ is empty and that
the set of attached prime ideals of any quotient of $A$ is contained
in $\Att_RA$. Also over a complete local ring $(R,\fm)$, by Matlis
Duality Theorem and using this description of attached primes of
$A$, it is straightforward to deduce the known fact that
$\Att_RA=\Ass_R(\Hom_R(A,E_R(R/\fm)))$. For basic theory concerning
attached prime ideals, we refer the reader to \cite[Section 6,
Appendix]{Mat}.

\begin{corollary} Let $(R,\fm)$ be a $d$-dimensional Cohen-Macaulay
local ring and $M$ and $N$ two finitely generated $R$-modules such that
either $\pd M$ or $\id N$ is finite. Then
$$\Att_R(H^d_{\fm}(M,N))=\Supp_RN\cap \Ass_RM.$$ In particular,
$H^d_{\fm}(M,N)=0$ if and only if $\Supp_RN\cap \Ass_RM=\emptyset$.
\end{corollary}

{\bf Proof.} Let $f:R\lo T$ be a ring homomorphism and $L$ a
$T$-module. Because for any $T$-module $K$, we have $\Ann_RK=
f^{-1}(\Ann_TK)$, it follows that $\Ass_RL=\{f^{-1}(\fp):\fp\in
\Ass_TL\}$. Also by the same reason, if $L$ is Artinian as an
$R$-module, then  $\Att_RL=\{f^{-1}(\fp):\fp\in \Att_TL\}$.
Moreover, one can easily check that the $\Supp_RN$ is precisely
consisting of the contractions to $R$ of the elements of
$\Supp_{\hat{R}}\hat{N}$. Hence, we may assume that $R$ is complete.

First assume that $\pd M<\infty$. Then, by Suzuki's Duality Theorem,
$H^d_{\fm}(M,N)\cong \Hom_R(N,M\otimes_R\omega_R)^{\vee}$. So, from
\cite[Exercise 1.2.27]{BH} and Lemma 3.3 ii), it follows that
$$\begin{array}{ll}
\Att_R(H^d_{\fm}(M,N))&=\Ass_R(\Hom_R(N,M\otimes_R\omega_R))\\
&=\Supp_RN\cap
\Ass_R(M\otimes_R\omega_R)\\
&=\Supp_RN\cap \Ass_RM.
\end{array}$$

Now, assume that $\id N<\infty$. Then, by the Herzog-Zamani Duality
Theorem, $$H_{\fm}^d(M,N)^{\vee}\cong
\Hom_R(\Hom_R(\omega_R,N),M).$$ Hence, from \cite[Exercise
1.2.27]{BH} and Lemma 3.3 iii), it follows that
$$\begin{array}{ll}
\Att_R(H^d_{\fm}(M,N))&=\Ass_R(\Hom_R(\Hom_R(\omega_R,N),M)))\\
&=\Supp_R(\Hom_R(\omega_R,N))\cap \Ass_RM\\
&=\Supp_RN\cap \Ass_RM. \  \  \Box
\end{array}
$$

Now, we are ready to present the first main result of this section.

\begin{theorem} Let $(R,\fm)$ be a Cohen-Macaulay local ring and
$M$ and $N$ two finitely generated $R$-modules such that either $\pd M$ or
$\id N$ is finite. Then $$\cd_{\fm}(M,N)=\dim R-\grade(\Ann_RN,M).$$
\end{theorem}

{\bf Proof.} Let $d:=\dim R$. First, assume that $\pd M<\infty$. Then by
Suzuki's Duality Theorem, $H^i_{\fm}(M,N)\cong
\Ext^{d-i}_{\hat{R}}(\hat{N},\hat{M}\otimes_{\hat{R}}\omega_{\hat{R}})
^{\vee}$. Hence by \cite[Proposition 1.2.10 b)]{BH} and Lemma 3.3
i), it follows that
$$\begin{array}{ll} \cd_{\fm}(M,N)&=\sup\{i:\Ext^{d-i}_{\hat{R}}(\hat{N},
\hat{M}\otimes_{\hat{R}}\omega_{\hat{R}})^{\vee}\neq 0\}\\
&=d-\inf\{j:\Ext^j_{\hat{R}}(\hat{N},\hat{M}\otimes_{\hat{R}}
\omega_{\hat{R}})\neq 0\}\\
&=d-\grade((\Ann_RN)\hat{R},\hat{M}\otimes_{\hat{R}}\omega_{\hat{R}})\\
&=d-\grade(\Ann_RN,M).
\end{array}$$

Now, assume that $\id N<\infty$. Then by the Herzog-Zamani Duality
Theorem  $$H^i_{\fm}(M,N)\cong
\Ext^{d-i}_{\hat{R}}(\Hom_{\hat{R}}(\omega_{\hat{R}},
\hat{N}),\hat{M}) ^{\vee}.$$ Therefore, in view of Lemma 3.3 iii),
the assertion follows by repeating the above argument. $\Box$

Next, in the case $R$ is a Cohen-Macaulay local ring, we improve the
dimension inequality $\dim N\leq \pd M+\dim (M\otimes_RN)$, which is
one of the famous consequence of the New Intersection Theorem.

\begin{corollary} Let $(R,\fm)$ be a Cohen-Macaulay local ring and
$M$ and $N$ two finitely generated $R$-modules. Assume that $\pd M$
is finite. Then $$\dim N\leq \dim R-\grade(\Ann_RN,M)\leq
\Gpd_NM+\dim (M\otimes_RN).$$
\end{corollary}

{\bf Proof.} By Theorem 3.5, $\cd_{\fm}(M,N)=\dim
R-\grade(\Ann_RN,M)$, and so $$\dim R-\grade(\Ann_RN,M)\leq
\Gpd_NM+\dim(M\otimes_RN),$$ by parts iv) and v) of Corollary 3.2.
On the other hand, since $\pd M$ is finite, any $M$-sequence is also
an $R$-sequence. Hence  $\grade(\Ann_RN,M)\leq \grade(\Ann_RN,R)=\Ht
(\Ann_RN)$, and so
$$\dim N=\dim R-\Ht (\Ann_RN)\leq \dim R-\grade(\Ann_RN,M). \Box$$

\begin{example} In Theorem 3.5, the Cohen-Macaulayness assumption on
$R$ is necessary. To see this, let $(R,\fm)$ be a non Cohen-Macaulay
local ring and $M, N$  finitely generated $R$-modules. Suppose that
$M$ is Cohen-Macaulay and $\pd M<\infty$. Recall that the height of
an ideal $\fa$ with respect to  $M$ is defined by $\Ht_M\fa:=\min\{ \dim
M_{\fp}:\fp\supseteq\fa\}.$ We have
$$\begin{array}{llll} \dim R-\grade(\Ann_RN,M)&>\depth R-\Ht_
M(\Ann_RN)\\
&=\pd M+\dim M-\Ht_M(\Ann_RN)\\
&=\pd M+\dim (M/(\Ann_RN)M)\\
&=\Gpd_NM+\dim (M\otimes_RN).
\end{array}$$
\end{example}

In the next result, we  compute the kernel of the epimorphism
$H_{\fm}^d(M,N)\lo H_{\fa}^d(M,N)$ in Corollary 3.2, when $d=\dim R$
and either $\pd M$ or $\id N$ is finite. Also, the following result
is crucial for proving our analogue of the Hartshorne-Lichtenbaum
Vanishing Theorem for generalized local cohomology modules.
In what follows, for an Artinian $R$-module $A$, we put
$<\fm>A:=\underset{i\in \mathbb{N}}\cap \fm^i A$.

\begin{theorem} Let $\fa$ be an ideal of a $d$-dimensional
Cohen-Macaulay local ring $(R,\fm)$ and $M$ and $N$ two finitely generated
$R$-modules. If either $\pd M$ or $\id N$ is finite, then there is a
natural isomorphism
$$H^d_{\fa}(M,N)\cong \frac{H^d_{\fm}(M,N)}{\underset{i\in
\mathbb{N}} \sum <\fm>(0:_{H^d_{\fm}(M,N)}\fa^n)}.$$
\end{theorem}

{\bf Proof.} We can assume that $R$ is complete. Let $\omega_R$ be
the canonical module of $R$ and
$(\cdot)^\vee:=\Hom_R(\cdot,E_R(R/\fm))$. Let $K$ and $L$ be
finitely generated $R$-modules such that $\pd L<\infty$ and set
$A:=K^\vee$. Suzuki's Duality Theorem asserts that
$$\Ext_R^d(K/\fa^n K, L\otimes_R\omega_R)\cong H_{\fm}^0
(L,K/\fa^n K)^\vee.$$ For a fixed integer $n\in \mathbb{N}$, choose
an integer $t:=t(n)\in \mathbb{N}$ such that
$$H_{\fm}^0(L,K/\fa^n K)\cong \Hom_R(R/\fm^t,\Hom_R(L,K/\fa^n K))\cong
\Hom_R(R/\fm^t\otimes_RL,K/\fa^n K)$$ and
$<\fm>(0:_A\fa^n)=\fm^t(0:_A\fa^n)$. It follows that
$$\begin{array}{lll} \Ext_R^d(K/\fa^n K,L\otimes_R\omega_R)
&\cong (R/\fm^t\otimes_RL)\otimes_R(K/\fa^n K)^\vee\\
&\cong L\otimes_R(R/\fm^t\otimes_R(0:_A\fa^n))\\
&\cong L\otimes_R(\frac{0:_A\fa^n}{<\fm>(0:_A\fa^n)}).
\end{array}$$
It is easy to check that
$\underset{n}{\varinjlim}\frac{0:_A\fa^n}{<\fm>(0:_A\fa^n)}=
\frac{A}{\underset{n\in \mathbb{N}}\sum <\fm>(0:_A\fa^n)}$. So
\cite[Lemma 3.1]{DS} yields that
$$  H^d_{\fa}(K,L\otimes_R\omega_R)\cong L\otimes_R
\frac{A}{\underset{n\in \mathbb{N}}\sum <\fm>(0:_A\fa^n)}\cong
\frac{L\otimes_RK^\vee} {\underset{n\in \mathbb{N}}\sum
<\fm>(0:_{L\otimes_RK^\vee}\fa^n)}. \ \  (*)$$ Assume that $\id
N<\infty$. Then by \cite[Theorem 4.3 ii)]{Su}, there is a natural
isomorphism $N\cong \Hom_R(\omega_R,N)\otimes_R\omega_R.$ On the
other hand,
$$\Hom_R(\omega_R,N)\otimes_RM^\vee\cong
\Hom_R(\Hom_R(\omega_R,N),M)^\vee,$$ and the later is isomorphic
with $H_{\fm}^d(M,N)$, by the Herzog-Zamani Duality Theorem. Since,
by \cite[Proposition 4.5 ii)]{Su}, $\pd
(\Hom_R(\omega_R,N))<\infty$, from $(*)$, it follows that
$$\begin{array}{lll} H^d_{\fa}(M,N)&\cong H^d_{\fa}(M,
\Hom_R(\omega_R,N)\otimes_R\omega_R)\\
&\cong \frac{H^d_{\fm}(M,N)}{\underset{i\in \mathbb{N}} \sum
<\fm>(0:_{H^d_{\fm}(M,N)}\fa^n)}.
\end{array}$$

Now, assume that $\pd M<\infty$. By \cite[Proposition 4.5 i)]{Su},
the functor $\cdot\otimes_R\omega_R$ is exact on the subcategory of
modules of finite projective dimension, and consequently it follows
that $\Tor^R_i(\omega_R,M)=0$ for all $i>0$. Hence, \cite[Lemma
2.5]{CD} and $(*)$ imply that
$$\begin{array}{lll} H^d_{\fa}(M,R)&\cong H^d_{\fa}(M
\otimes_R\omega_R,\omega_R)\\
&\cong \frac{(M\otimes_R\omega_R)^\vee}{\underset{i\in \mathbb{N}}
\sum
<\fm>(0:_{(M\otimes_R\omega_R)^\vee}\fa^n)}\\
&\cong \frac{H^d_{\fm}(M,R)}{\underset{i\in \mathbb{N}} \sum
<\fm>(0:_{H^d_{\fm}(M,R)}\fa^n)}.
\end{array}$$
Note that by Suzuki's Duality Theorem
$(M\otimes_R\omega_R)^\vee\cong H^d_{\fm}(M,R)$. On the other hand,
by Corollary 3.2 v), the functors $H^d_{\fa}(M,\cdot)$ and
$H^d_{\fm}(M,\cdot)$  are right exact, and so using \cite[Lemma
3.1]{DS} once more, yields that
$$
\begin{array}{lll} H^d_{\fa}(M,N)&\cong
N\otimes_R\frac{H^d_{\fm}(M,R)}{\underset{i\in \mathbb{N}} \sum
<\fm>(0:_{H^d_{\fm}(M,R)}\fa^n)}\\
&\cong \frac{H^d_{\fm}(M,N)}{\underset{i\in \mathbb{N}} \sum
<\fm>(0:_{H^d_{\fm}(M,N)}\fa^n)}. \ \  \Box
\end{array}
$$

\begin{corollary} Let $\fa$ be an ideal of a $d$-dimensional
Cohen-Macaulay local ring $(R,\fm)$ and $M$ and $N$ two finitely generated
$R$-modules. If either $\pd M$ or $\id N$ is finite, then
$$\Att_{\hat{R}}(H^d_{\fa}(M,N))=\{\fp\in \Supp_{\hat{R}}\hat{N}\cap
\Ass_{\hat{R}}\hat{M}:\dim (\hat{R}/\fa \hat{R}+\fp)=0\}.$$
\end{corollary}

{\bf Proof.} We can assume that $R$ is complete. Let
$H^d_{\fm}(M,N)=A_1+ \dots +A_n$ be a minimal secondary
representation of $H^d_{\fm}(M,N)$. We order the set
$\Att_R(H^d_{\fm}(M,N))=\{\fp_1,\dots , \fp_n\}$ such that for an
integer $0\leq l\leq n$, $\dim (R/\fa+\fp_i)>0$ for all $1\leq i\leq
l$, while $\dim (R/\fa+\fp_i)=0$ for all $l+1\leq i\leq n$. Then by
\cite[Theorem 2.8]{DS}, $A_1+\dots +A_l$ is a minimal secondary
representation of $B:=\sum_{n\in
\mathbb{N}}<\fm>(0:_{H^d_{\fm}(M,N)}\fa^n)$.  Now, it is easy to see
that $\Sigma_{i=l+1}^n(A_i+B)/B$ is a minimal secondary
representation of $H^d_{\fa}(M,N)\cong H^d_{\fm}(M,N)/B$. But
$(A_i+B)/B\cong A_i/(A_i\cap B)$ is $\fp_i$-secondary for all
$l+1\leq i\leq n$. This finishes the proof by Corollary 3.4. $\Box$

When $R$ is a Gorenstein local ring, an analogue of the
Hartshorne-Lichtenbaum Vanishing Theorem for generalized local
cohomology modules was  established in \cite[Lemma 2.5]{DST}. Now,
we weaken the assumption on $R$ to the Cohen-Macaulayness.

\begin{corollary} Let $\fa$ be an ideal of a $d$-dimensional
Cohen-Macaulay local ring $(R,\fm)$ and $M$ and $N$ two finitely generated
$R$-modules. Assume that either $\pd M$ or $\id N$ is finite. Then
the following are equivalent:
\begin{enumerate}
\item[i)] $H^d_{\fa}(M,N)=0$.
\item[ii)] $H^d_{\fm}(M,N)=\displaystyle{\sum_{n\in
\mathbb{N}}}<\fm>(0:_{H^d_{\fm}(M,N)}\fa^n)$.
\item[iii)] For any integer $l\in \mathbb{N}$, there exists an
$n=n(l)\in \mathbb{N}$ such that
$$0:_{H^d_{\fm}(M,N)}\fa^l\subseteq
<\fm>(0:_{H^d_{\fm}(M,N)}\fa^n).$$
\item[iv)] $\dim \hat{R}/\fa \hat{R}+\fp>0$ for all
$\fp\in \Supp_{ \hat{R}}\hat{N}\cap \Ass_{\hat{R}}\hat{M}$.
\end{enumerate}
\end{corollary}

{\bf Proof.} i), ii) and iv) are equivalent by Theorem 3.8 and
Corollary 3.9, while the equivalence of ii) and iii) follows by
\cite[Corollary 2.5]{DS}. $\Box$

%%%%%%%%%%%%%%%%%%%%%%%%%%%%%%%%%%%%%%%%%%%%%%%%%%%%%%%%%%%%%%%%%

\end{document}